\pdfoutput=1
\documentclass{amsart}
\usepackage[utf8]{inputenc}

\usepackage[margin=1in]{geometry}
\usepackage[expansion=false]{microtype}

\usepackage{amsmath, amsfonts, amssymb, amsthm, amsopn}
\usepackage{eucal}

\usepackage{tikz}
\usetikzlibrary{diagrams}

\usepackage[pdfusetitle,unicode,hidelinks]{hyperref}

\newcommand{\ie}{{\sfcode`\.1000 i.e.}}

\numberwithin{equation}{section}

\theoremstyle{plain}
\newtheorem{theorem}[equation]{Theorem}
\newtheorem{proposition}[equation]{Proposition}

\newtheorem{lemma}[equation]{Lemma}
\newtheorem{corollary}[equation]{Corollary}

\theoremstyle{definition}
\newtheorem{definition}[equation]{Definition}
\newtheorem{example}[equation]{Example}
\newtheorem{remark}[equation]{Remark}

\let\scr=\mathcal
\let\bb=\mathbf

\let\phi=\varphi

\def\ph{\mathord-}

\def\Z{\bb Z}
\def\Q{\bb Q}

\def\1{\bb 1}

\def\Set{\scr S\mathrm{et}}

\def\et{{\mathrm{\acute et}}}

\def\Zar{{\textrm{Zar}}}

\def\fc{\mathrm{fc}}

\def\op{\mathrm{op}}

\let\into=\hookrightarrow

\let\comma=/

\def\id{\mathrm{id}}

\def\Tors{\mathrm{Tors}}

\def\B{\beta}
\def\T{\theta}

\def\Top{\scr T\mathrm{op}}

\def\Cat{\scr C\mathrm{at}}

\def\Ab{\scr A\mathrm{b}}
\def\ET{{\scr E}\mathrm{t}}

\def\HC{\mathrm{HC}}
\def\EM{\mathcal{EM}}
\def\h{\mathrm{h}}

\DeclareMathOperator{\Hom}{Hom}
\DeclareMathOperator{\Map}{Map}
\DeclareMathOperator{\Fun}{Fun}

\DeclareMathOperator{\Pro}{Pro}
\DeclareMathOperator{\Ind}{Ind}

\DeclareMathOperator{\Shv}{Shv}

\DeclareMathOperator{\Ex}{Ex}

\let\lim=\relax
\DeclareMathOperator*{\lim}{lim}

\DeclareMathOperator*{\colim}{colim}
\DeclareMathOperator*{\hocolim}{hocolim}

\title{Higher Galois theory}
\author{Marc Hoyois}
\date{\today}
\address{Department of Mathematics, Massachusetts Institute of Technology, Cambridge, MA, USA}
\email{hoyois@mit.edu}
\urladdr{\url{http://math.mit.edu/~hoyois/}}

\begin{document}

\maketitle

\begin{abstract}
We generalize toposic Galois theory to higher topoi.
We show that locally constant sheaves in a locally $(n-1)$-connected $n$-topos are equivalent to representations of its fundamental pro-$n$-groupoid, and that the latter can be described in terms of Galois torsors.
We also show that finite locally constant sheaves in an arbitrary $\infty$-topos are equivalent to finite representations of its fundamental pro-$\infty$-groupoid. Finally, we relate the fundamental pro-$\infty$-groupoid of $1$-topoi to the construction of Artin and Mazur and, in the case of the étale topos of a scheme, to its refinement by Friedlander.
\end{abstract}

\tableofcontents

\section{Introduction}

The goal of this paper is to generalize Galois theory, as it is understood in the context of Grothendieck topoi, to higher topoi.
We assume that the reader is familiar with the theory of $\infty$-topoi as it is developed in \cite{HTT}.

In \S\ref{sec:shape}, we review the notion of shape, or fundamental pro-$\infty$-groupoid, of an $\infty$-topos. For the original accounts, we refer the reader to Toën and Vezzosi \cite[\S5.3]{TV} and to Lurie \cite[\S7.1.6]{HTT}.

In \S\ref{sec:galois}, we prove our generalization of (infinite) Galois theory to $n$-topoi, for $0\leq n\leq\infty$. Specializing to $n=1$ recovers classical results of Moerdijk \cite{Moerdijk} and of Kennison \cite{Kennison}, but our proofs are quite different as they make essential use of higher topos theory even in that case (specifically, of the theory of stacks in groupoids).
The main result is that locally constant sheaves in a locally $(n-1)$-connected $n$-topos are equivalent to representations of its fundamental pro-$n$-groupoid. For locally $n$-connected (and $1$-localic) $n$-topoi, this statement already appears in Grothendieck's famous letter to Breen \cite[(C)]{GrothendieckBreen},
but most of our work deals with the weakening of the local connectedness assumption, which is necessary to recover the Galois theory of fields. As consequences, we obtain a model for the fundamental pro-$n$-groupoid as a cofiltered diagram indexed by Galois torsors, and we deduce the following recognition principle: a locally $(n-1)$-connected $n$-topos that is generated by its locally constant sheaves is equivalent to the classifying $n$-topos of a pro-$n$-groupoid with $(n-1)$-connected transition maps.

In \S\ref{sec:finitegalois}, we generalize finite Galois theory to arbitrary $\infty$-topoi: we show that finite locally constant sheaves in an $\infty$-topos are equivalent to finite representations of its fundamental pro-$\infty$-groupoid.

In \S\ref{sec:examples}, we observe that if $\scr C$ is a Grothendieck site with a final object, then a coherent refinement of the construction of Artin and Mazur \cite{AM} is a model for the shape of the $\infty$-topos of hypercomplete sheaves on $\scr C$. Specializing to the étale site of a locally connected scheme $X$, we show that the étale topological type of $X$ defined by Friedlander \cite{Friedlander:1982} is a model for the shape of the $\infty$-topos of hypercomplete étale sheaves on $X$.

\subsubsection*{Acknowledgments.}
I am grateful to Tim Porter for bringing Grothendieck's letter \cite{GrothendieckBreen} to my attention, and to Eduardo Dubuc for interesting discussions about his construction of the fundamental pro-groupoid of a topos in \cite{Dubuc,Dubuc2}.

\subsubsection*{Notation.}
\leavevmode

\begin{tabular}[t]{ll}
	$\infty$ & $(\infty,1)$ \\
	$\scr S$      & $\infty$-category of small $\infty$-groupoids \\
	$\Cat$ & $\infty$-category of (possibly large) $\infty$-categories \\
	$\Top$  & $\infty$-category of $\infty$-topoi and geometric morphisms \\
	$\1$ & final object in an $\infty$-topos \\
	$\Fun(\scr C,\scr D)$ & $\infty$-category of functors from $\scr C$ to $\scr D$ \\
	$\Map(X,Y)$ & $\infty$-groupoid of maps from $X$ to $Y$ in an $\infty$-category \\
	$\scr C_{\comma X}$, $\scr C_{X\comma}$ & overcategory, undercategory \\
	$\scr C_{\leq n}$ & subcategory of $n$-truncated objects \\
	$\h\scr C$ & homotopy category \\
	$\Set_\Delta$ & category of simplicial sets \\
	$\Map_{\Delta}(X,Y)$ & simplicial set of maps from $X$ to $Y$ in a $\Set_\Delta$-enriched category \\
\end{tabular}

\section{Preliminaries on shapes}
\label{sec:shape}

The shape of an $\infty$-topos is a pro-object in the $\infty$-category $\scr S$ of $\infty$-groupoids. We begin with a brief review of pro-objects in $\infty$-categories.

Let $\scr C$ be an $\infty$-category. The $\infty$-category $\Pro(\scr C)$ of pro-objects in $\scr C$ and the Yoneda embedding $j\colon\scr C\to\Pro(\scr C)$ are defined by the following universal property:
\begin{enumerate}
	\item The $\infty$-category $\Pro(\scr C)$ admits small cofiltered limits;
	\item Let $\scr D$ be an $\infty$-category which admits small cofiltered limits, and let $\Fun'(\Pro(\scr C),\scr D)$ be the full subcategory of $\Fun(\Pro(\scr C),\scr D)$ spanned by the functors that preserve small cofiltered limits. Then $j$ induces an equivalence $\Fun'(\Pro(\scr C),\scr D)\simeq\Fun(\scr C,\scr D)$.
\end{enumerate}

The existence of $\Pro(\scr C)$ in general is a special case of \cite[Proposition 5.3.6.2]{HTT}.
If $\scr C$ is accessible and admits finite limits, then $\Pro(\scr C)$ can be identified with the full subcategory of $\Fun(\scr C,\scr S)^\op$ spanned by the left exact \emph{accessible} functors \cite[Proposition 3.1.6]{DAG13}. 
The universal property is then explicitly realized as follows. Recall that $\scr S$ is the base of the universal left fibration $u\colon\scr S_{\ast\comma}\to\scr S$. Thus, any functor $f\colon \scr C\to\scr S$ classifies a left fibration $\scr C_{f\comma}\to\scr C$ given by the cartesian square
\begin{tikzmath}
	\diagram{\scr C_{f\comma} & \scr S_{\ast\comma} \\ \scr C & \scr S\rlap.\\};
	\arrows (11-) edge (-12) (11) edge (21) (12) edge node[right]{$u$} (22) (21-) edge node[below]{$f$} (-22);
\end{tikzmath}
The condition that $f$ is left exact implies that $\scr C_{f\comma}$ is cofiltered \cite[Remark 5.3.2.11]{HTT}, and the condition that $f$ is accessible implies that $\scr C_{f\comma}$ is accessible \cite[Proposition 5.4.6.6]{HTT}. Accessible $\infty$-categories have small left cofinal subcategories since for $\scr A$ small, $\scr A\subset\Ind_\kappa(\scr A)$ is left cofinal.\footnote{We follow the convention of Bousfield and Kan \cite[XI \S9]{BK} regarding left vs.\ right cofinality. Unfortunately, the opposite convention is used in \cite{HA}.} Thus, diagrams indexed by $\scr C_{f\comma}$ will have a limit in any $\infty$-category $\scr D$ that admits small cofiltered limits. In this way any functor $\scr C\to\scr D$ lifts to a functor $\Pro(\scr C)\to\scr D$, sending $f$ to the limit of the composition $\scr C_{f\comma}\to \scr C\to\scr D$. Note that the cofiltered diagram $u\colon \scr C_{f\comma}\to\scr C$ \emph{corepresents} the pro-object $f$, in the sense that
\[f(K)\simeq\colim_{X\in\scr C_{f\comma}}\Map(uX,K).\]
 
\begin{remark}
	If $\scr C$ is an $n$-category, the $\infty$-category $\Pro(\scr C)$ is again an $n$-category \cite[Remark 5.3.5.6]{HTT}. If $\scr C$ is an ordinary category (viewed as an $\infty$-category with discrete mapping spaces), it follows that $\Pro(\scr C)$ is equivalent to the category of pro-objects defined in \cite[Exposé I, \S8.10]{SGA4-1}, since both satisfy the same universal property \cite[Exposé I, Proposition 8.7.3]{SGA4-1}.
\end{remark}

\begin{remark}
	By \cite[Proposition 5.3.1.16]{HTT}, any pro-object in an $\infty$-category $\scr C$ can be corepresented by a diagram $\scr I\to\scr C$ where $\scr I$ is a small cofiltered \emph{poset}. Using this fact, one can show that, if $\scr C$ is the underlying $\infty$-category of a proper model category $\scr M$, then $\Pro(\scr C)$ is the underlying $\infty$-category of the \emph{strict} model structure on $\Pro(\scr M)$ defined in \cite{Isaksen:2007}; we will not need this model in the sequel.
\end{remark}

Let $G\colon\scr D\to\scr C$ be an accessible functor between presentable $\infty$-categories, and let $F\colon\scr C\to\Fun(\scr D,\scr S)^\op$ be the ``formal'' left adjoint to $G$. By the adjoint functor theorem, $F$ factors through $\scr D$ if and only if $G$ preserves small limits. Clearly, $F$ factors through $\Pro(\scr D)$ if and only if $G$ preserves finite limits. In this case the functor $F\colon\scr C\to\Pro(\scr D)$ is called the \emph{pro-left adjoint} to $G$. Its extension $\Pro(\scr C)\to\Pro(\scr D)$ is a genuine left adjoint to $\Pro(G)$.

If $f\colon\scr Y\to\scr X$ is a geometric morphism of $\infty$-topoi, we will write $f_\ast\colon\scr Y\to\scr X$ for the direct image functor and $f^\ast\colon\scr X\to\scr Y$ for the left adjoint of the latter. Since $f^\ast$ preserves finite limits, it admits a pro-left adjoint $f_!\colon\scr Y\to\Pro(\scr X)$ given by
\[f_!(Y)(X)\simeq\Map_{\scr Y}(Y, f^\ast X).\]

If $\scr X$ is an $\infty$-topos, we will usually denote by $\pi\colon\scr X\to\scr S$ the unique geometric morphism to $\scr S$, given informally by $\pi_\ast(X)=\Map(\1,X)$. 
The following definition is due to Toën and Vezzosi \cite[\S5.3]{TV} and to Lurie \cite[\S7.1.6]{HTT}, following ideas of Grothendieck \cite{GrothendieckBreen}.

\begin{definition}
	Let $\scr X$ be an $\infty$-topos and $\pi\colon\scr X\to\scr S$ the unique geometric morphism.
	The \emph{fundamental pro-$\infty$-groupoid} or \emph{shape} of $\scr X$, denoted by $\Pi_\infty\scr X$, is the pro-$\infty$-groupoid defined by
	\[\Pi_\infty\scr X=\pi_!\1\in\Pro(\scr S).\]
\end{definition}

As a left exact functor $\scr S\to\scr S$, $\Pi_\infty\scr X$ is the composition $\pi_\ast\pi^\ast$, \ie, it sends an $\infty$-groupoid to the global sections of the associated constant sheaf. Note that $\pi_! X\simeq\Pi_\infty(\scr X_{\comma X})$ since $X\simeq\rho_!\1$ where $\rho\colon\scr X_{\comma X}\to \scr X$ is the canonical geometric morphism. In other words, the functor $\Pi_\infty\colon\Top\to\Pro(\scr S)$ simultaneously extends the functors $\pi_!\colon\scr X\to\Pro(\scr S)$ for all $\infty$-topoi $\scr X$.

\begin{example}
	\label{ex:CW}
	If $T$ is a topological space homotopy equivalent to a CW complex, the shape of $\Shv(T)$ is the weak homotopy type of $T$ \cite[Remarks A.1.4 and A.4.7]{HA}. As we will see in Proposition~\ref{prop:observation} below, this is a refinement of the fact that the singular cohomology of $T$ with coefficients in a local system coincides with its sheaf cohomology with coefficients in the corresponding locally constant sheaf.
\end{example}

\begin{definition}
	Let $\scr X$ be an $\infty$-topos.
	A \emph{torsor} $(A,\chi)$ in $\scr X$ is an $\infty$-groupoid $A\in\scr S$ together with a map $\chi\colon \1\to\pi^*A$ in $\scr X$. We denote by
\[\Tors(\scr X)=\scr X_{\1\comma}\times_{\scr X}\scr S\]
the $\infty$-category of torsors in $\scr X$ and by
\[\Tors(\scr X,A)=\Tors(\scr X)\times_{\scr S}\{A\}\simeq\Map_{\scr X}(\1,\pi^*A)\]
the $\infty$-groupoid of \emph{$A$-torsors}.	
\end{definition}

  By descent, an $A$-torsor is equivalently an action $P\colon A\to\scr X$ of the $\infty$-groupoid $A$ in $\scr X$ which is \emph{principal} in the sense that $\colim_{\alpha\in A} P_\alpha\simeq\1$. 
  Note that the projection $\Tors(\scr X)\to\scr S$ is the left fibration classified by $\pi_*\pi^*\colon\scr S\to\scr S$. In particular, it corepresents the shape of $\scr X$.

\begin{proposition}
	Let $f\colon\scr Y\to\scr X$ be a geometric morphism of $\infty$-topoi. The following conditions are equivalent:
	\begin{enumerate}
		\item $f$ is a shape equivalence, \ie, $\Pi_\infty(f)$ is an equivalence in $\Pro(\scr S)$.
		\item $f^*\colon\Tors(\scr X)\to\Tors(\scr Y)$ is an equivalence of $\infty$-categories.
		\item For every $A\in\scr S$, $f^*\colon\Tors(\scr X,A)\to\Tors(\scr Y,A)$ is an equivalence of $\infty$-groupoids.
	\end{enumerate}
\end{proposition}

\begin{proof}
	As we recalled above, the Grothendieck construction provides an equivalence between accessible functors $\scr S\to\scr S$ and accessible left fibrations over $\scr S$, under which $\Pi_\infty\scr X$ corresponds to $\Tors(\scr X)$. Explicitly, we have natural equivalences
	\[\Pi_\infty\scr X\simeq \lim_{(A,\chi)\in\Tors(\scr X)}j(A)\quad\text{and}\quad \Tors(\scr X)\simeq\scr S_{\Pi_\infty\scr X\comma}.\]
This proves $(1)\Leftrightarrow (2)$. A morphism of left fibrations is an equivalence if and only if it is a fiberwise equivalence, so $(2)\Leftrightarrow(3)$.
\end{proof}

We will now construct a right adjoint to the functor $\Pi_\infty$.
Recall that, for any $\infty$-topos $\scr X$, there is a fully faithful functor
\[
\scr X\into \Top_{\comma \scr X},\quad U\mapsto \scr X_{\comma U}.
\]
When $\scr X=\scr S$, we denote this functor by
\[
\B\colon \scr S\into \Top,
\]
and we call $\B A=\scr S_{\comma A}\simeq\Fun(A,\scr S)$ the \emph{classifying $\infty$-topos} of the $\infty$-groupoid $A$.
The $\infty$-topos $\B A$ classifies $A$-torsors, meaning that there is an equivalence of $\infty$-categories
\begin{equation}\label{eqn:Xtorsors}
	\Fun_{\Top}(\scr X,\B A)\simeq \Tors(\scr X,A),
\end{equation}
for every $\infty$-topos $\scr X$ (this is a special case of \cite[Corollary 6.3.5.6]{HTT}). The diagonal map $A\to A\times A$ in $\scr S_{\comma A}$ is the universal $A$-torsor.

Since the $\infty$-category $\Top$ admits small cofiltered limits, $\B$ extends to a functor
\[\B\colon\Pro(\scr S)\to\Top.\]
Explicitly, let $X\in\Pro(\scr S)$ be a pro-$\infty$-groupoid given in the form of a cofiltered diagram $X\colon\scr I\to\scr S$. The classifying $\infty$-topos $\B(X)$ is then the limit of the cofiltered diagram
\[\B\circ X\colon\scr I\to\Top.\]
Recall that limits of cofiltered diagrams in $\Top$ are created by the forgetful functor $\Top\to\Cat$ \cite[Theorem 6.3.3.1]{HTT}. Thus, an object $L\in\B(X)$ is a family of objects $L_i\in\scr \B(X_i)$  together with coherent equivalences $f_*L_i\simeq L_j$ for all arrows $f\colon i\to j$ in $\scr I$.

\begin{remark}
	While the functor $\B$ is fully faithful and colimit-preserving on $\scr S$, its extension to $\Pro(\scr S)$ does not have either property: if $X$ is a pro-set whose limit is empty, it is clear that $\B(X)$ is the empty $\infty$-topos, and if $\B$ preserved pushouts of pro-$\infty$-groupoids, it would imply that cofiltered limits preserve pushouts in $\scr S$. Nevertheless, one can show that $\B$ preserves arbitrary coproducts of pro-$\infty$-groupoids (and, more generally, colimits of ``levelwise'' diagrams indexed by $\infty$-groupoids).
\end{remark}

\begin{remark}
	For $X\in\Pro(\scr S)$, the $\infty$-topos $\B X$ is typically not hypercomplete. For example, it is shown in \cite[Warning 7.2.2.31]{HTT} that, for $p$ prime, the classifying $\infty$-topos $\B(B\Z_p)$ is not hypercomplete.
\end{remark}

If $X\colon\scr I\to\scr S$ is a pro-$\infty$-groupoid and $\scr C$ is an $\infty$-category, it is convenient to write
\[\Fun(X,\scr C)=\colim_{i\in\scr I}\Fun(X_i,\scr C)\]
---in other words, we view both $X$ and $\scr C$ as pro-$\infty$-categories.
The inverse image functors $\B(X_i)\to\B(X)$ then induce a canonical functor
\[\Fun(X,\scr S)\to\B(X);\]
an object $L\in\B(X)$ is called \emph{split} if it belongs to its essential image, and if $L$ is the image of $K\in\B(X_i)$ we say that $L$ is \emph{represented} by $K$. An arbitrary $L\in\B(X)$ is the filtered colimit of the split objects represented by $L_i$ for $i\in\scr I$ \cite[Lemma 6.3.3.6]{HTT}.

\begin{definition}
	If $\scr X$ is an $\infty$-topos, $\B(\Pi_\infty\scr X)$ is the \emph{$\infty$-topos of local systems on $\scr X$}.
\end{definition}

Since $\Pi_\infty\scr X$ is corepresented by the forgetful functor $\Tors(\scr X)\to\scr S$, a local system on $\scr X$ is a family of objects $L_{(X,x)}\in\B(X)$ indexed by pairs $(X,x)$ where $X\in\scr S$ and $x\colon\1\to\pi^\ast X$ is an $X$-torsor in $\scr X$, together with coherent equivalences $f_\ast L_{(X,x)}\simeq L_{(Y,y)}$ for all morphisms of torsors $f\colon(X,x)\to (Y,y)$.

It follows from~\eqref{eqn:Xtorsors} that we have an adjunction
\begin{tikzmath}
	\diagram{\Top & \Pro(\scr S)\rlap;\\};
	\arrows (11-) edge[vshift=\dbl] node[above=\dbl]{$\Pi_\infty$} (-12) (-12) edge[vshift=\dbl] node[below=\dbl]{$\B$} (11-);
\end{tikzmath}
we will denote by $\phi\colon\scr X\to\B(\Pi_\infty\scr X)$ its unit. If $L$ is a local system on $\scr X$, the object $\phi^\ast L\in\scr X$ will be called the \emph{underlying sheaf}  of $L$. We wish to describe $\phi^\ast L$ more explicitly. Since a general local system is a filtered colimit of split local systems, it suffices to describe the underlying sheaves of the latter. Let $(X,x\colon \1\to\pi^\ast X)\in\Tors(\scr X)$ be a torsor. This determines a geometric morphism $f\colon\scr X\to\B(X)$. If $L$ is represented by $K\in\B(X)$, we therefore have $\phi^\ast L\simeq f^\ast K$. This means that $\phi^\ast L$ fits in a cartesian square
\begin{tikzmath}
	\diagram{\phi^\ast L & \pi^\ast K \\ \1 & \pi^\ast X\rlap.\\};
	\arrows (11-) edge (-12) (11) edge (21) (12) edge (22) (21-) edge node[above]{$x$} (-22);
\end{tikzmath}

In addition to preserving colimits (being left adjoint), the functor $\Pi_\infty\colon \Top\to\Pro(\scr S)$ also preserves some interesting limits:

\begin{proposition}
	\label{prop:lim}
	\leavevmode
	\begin{enumerate}
	\item If $\scr X$ is the limit of a cofiltered diagram of proper $\infty$-topoi $(\scr X_i)$ with proper transition morphisms, then $\Pi_\infty\scr X\simeq\lim_i\Pi_\infty\scr X_i$.
	\item If $\scr X$ and $\scr Y$ are proper $\infty$-topoi, then $\Pi_\infty(\scr X\times\scr Y)\simeq \Pi_\infty\scr X\times\Pi_\infty\scr Y$.
	\end{enumerate}
\end{proposition}

\begin{proof}
	(1) Let $\pi\colon\scr X\to \scr S$ and $\pi_i\colon\scr X_i\to\scr S$ be the unique geometric morphisms.
	Since proper geometric morphisms preserve filtered colimits \cite[Remark 7.3.1.5]{HTT}, the canonical map
	\[
	\colim_i \pi_{i*}\pi_i^*\to \pi_*\pi^*
	\]
	is an equivalence. This exactly says that $\Pi_\infty\scr X\simeq \lim_i\Pi_\infty\scr X_i$. 
	
	(2) The properness of $\scr Y$ implies, by proper base change, that
	\[
	\Pi_\infty(\scr X\times\scr Y)\simeq \Pi_\infty\scr X\circ \Pi_\infty\scr Y.
	\]
	The properness of $\scr X$ implies that $\Pi_\infty\scr X$ preserves filtered colimits.
	To conclude, note that if $X,Y\in\Pro(\scr S)$ and $X\colon\scr S\to\scr S$ preserves filtered colimits, then $X\circ Y\simeq X\times Y$.
\end{proof}

\begin{example}\label{ex:CHT}
	If a topological space $T$ is the limit of a cofiltered diagram $(T_i)$ of compact Hausdorff spaces, then
	\[
	\Pi_\infty\Shv(T)\simeq \lim_i\Pi_\infty \Shv(T_i).
	\]
	This is a consequence of Proposition~\ref{prop:lim} (1) and the following facts: passing to locales preserves the limit of the diagram $(T_i)$, and maps between compact Hausdorff spaces induce proper morphisms of $\infty$-topoi \cite[Theorem 7.3.1.16]{HTT}.
	Together with Example~\ref{ex:CW}, we deduce for example that, if $T$ is a compact subset of the Hilbert cube, then the shape of $\Shv(T)$ is corepresented by the diagram $\mathbb R_{>0}\to \scr S$ sending $\epsilon$ to the weak homotopy type of the $\epsilon$-neighborhood of $T$.
\end{example}

\begin{remark}
	Since any compact Hausdorff space is a cofiltered limit of compact metrizable ANRs, which have the homotopy type of CW complexes, Examples \ref{ex:CW} and~\ref{ex:CHT} show that the shape theory of $\infty$-topoi completely subsumes the classical shape theory of compact Hausdorff spaces \cite{MS}.
\end{remark}

For $n\geq -2$, consider the adjunction
\begin{tikzmath}
	\diagram{\Pro(\scr S) & \Pro(\scr S_{\leq n})\rlap.\\};
	\arrows (11-) edge[vshift=\dbl] node[above]{$\tau_{\leq {n}}$} (-12) (-12) edge[c->,vshift=\dbl] (11-);
\end{tikzmath}
We denote by $\Pi_n\colon\Top\to\Pro(\scr S_{\leq n})$ the composition $\tau_{\leq n}\circ\Pi_\infty$.
Note that $\Pi_n\scr X$ is corepresented by the forgetful functor $\Tors_{\leq n}(\scr X)\to\scr S_{\leq n}$, where $\Tors_{\leq n}(\scr X)=\Tors(\scr X)\times_{\scr S}\scr S_{\leq n}$.
The description of torsors in terms of principal actions shows that the functor $\Pi_n$ factors through the reflective subcategory $\Top_n\subset\Top$ of $n$-localic $\infty$-topoi. Moreover, if $X\in\Pro(\scr S_{\leq n})\subset\Pro(\scr S)$, then $\B X$ is a limit of $n$-localic $\infty$-topoi and hence is $n$-localic. Thus, we obtain an induced adjunction
\[
\Pi_n :\Top_n\rightleftarrows \Pro(\scr S_{\leq n}):\B
\]
between $n$-topoi and pro-$n$-groupoids.

\begin{remark}
	If $\scr X$ is a $1$-topos, $\Pi_1\scr X$ is equivalent to the fundamental pro-groupoid defined by Dubuc in \cite[\S8]{Dubuc2}, since both classify torsors.
	 Earlier constructions of $\Pi_1\scr X$ were given in special cases, notably by Kennison \cite{Kennison} (for $\scr X$ connected), by Bunge \cite{Bunge} (for $\scr X$ locally connected), and by Barr and Diaconescu \cite{BarrDiaconescu} (for $\scr X$ locally simply connected).
\end{remark}

By definition, $\Pi_\infty\scr X$ computes the cohomology of $\scr X$ with constant coefficients. Our next observation is that this can be extended to cohomology with coefficients in the underlying sheaf of a split local system.
In \S\ref{sec:galois} we will show that there are many such sheaves: for example, if $\scr X$ is locally connected, $A$ is a locally constant sheaf of abelian groups, and $n\geq 0$, then the Eilenberg–Mac Lane sheaf $K(A,n)$ is the underlying sheaf of a split local system.

\begin{proposition}\label{prop:observation}
	Let $\scr X$ be an $\infty$-topos and let $L\in\Fun(\Pi_\infty\scr X,\scr S)$. 
	Denote by $\tilde L$ the image of $L$ in $\B(\Pi_\infty\scr X)$.
	Then $\phi^*$ induces an equivalence
	\[\Map_{\Fun(\Pi_\infty\scr X,\scr S)}(\ast,L)\simeq\Map_{\scr X}(\1,\phi^*\tilde L).\]
\end{proposition}

\begin{proof}
	Suppose that $L$ comes from the object $K\in\scr S_{\comma X}$ labeled by the torsor $(X,x)$ in $\scr X$. The proposition follows by comparing the two cartesian squares
	\[
	\begin{tikzpicture}[baseline=(24.base)]
		\diagram{\Map_{\Fun(\Pi_\infty\scr X,\scr S)}(\ast,L) & \Map(\Pi_\infty\scr X,K) & \pi_\ast\phi^\ast\tilde L & \pi_\ast\pi^\ast K \\ \ast & \Map(\Pi_\infty\scr X, X)\rlap, & \ast & \pi_\ast\pi^\ast X\rlap, \\};
		\arrows (11-) edge (-12) (11) edge (21) (21-) edge (-22) (12) edge (22)
		(13-) edge (-14) (13) edge (23) (23-) edge (-24) (14) edge (24);
	\end{tikzpicture}
	\]
	in which the lower horizontal maps are induced by $x$.
\end{proof}

\section{Galois theory}
\label{sec:galois}

Classical Galois theory states that the étale topos $\scr X$ of a field $k$ is equivalent to the classifying topos of the absolute Galois group of $k$.
More precisely:
\begin{enumerate}
	\item For any separable closure $k^s$ of $k$, there is a canonical equivalence of pro-groupoids
	\[B\mathrm{Gal}(k^s/k)\simeq \Pi_1\scr X.\]
	\item The geometric morphism $\phi\colon\scr X\to\B(\Pi_1\scr X)$ is an equivalence of topoi and identifies locally constant sheaves with split local systems.
\end{enumerate}
Statement (2) is true more generally for any locally connected topos $\scr X$ generated by its locally constant object. There is also an analog of statement (1) for any locally connected topos.
In this section we prove the $n$-toposic generalizations of these results, for $0\leq n\leq\infty$.
The case $n=\infty$ is treated in \cite[\S{}A.1]{HA}, but the case of finite $n$ is more complicated.

\begin{definition}
	Let $\scr X$ be an $\infty$-topos. An object $X\in\scr X$ is called \emph{locally constant} if there exists an effective epimorphism $\coprod_\alpha U_\alpha\to\1$ such that $X$ is constant over each $U_\alpha$, \ie, such that $X\times U_\alpha\simeq \pi^\ast X_\alpha\times U_\alpha$ for some $X_\alpha\in\scr S$.
\end{definition}

\begin{definition}
Let $-2\leq n\leq \infty$.
A geometric morphism $f\colon \scr Y\to\scr X$ is called \emph{$n$-connected} if $f^*$ is fully faithful on $n$-truncated objects.
An $\infty$-topos $\scr X$ is called \emph{locally $n$-connected} if $\pi^*\colon \scr S_{\leq n}\to \scr X_{\leq n}$ preserves infinite products, or, equivalently, if its pro-left adjoint is a genuine left adjoint.	
\end{definition}

 Note that every $\infty$-topos is locally $(-1)$-connected, since $\scr S_{\leq -1}=\{\emptyset\to\ast\}$.

\begin{proposition}\label{prop:lconst}
	Let $\scr X$ be a locally $\infty$-connected $\infty$-topos. Then $\phi\colon\scr X\to\B(\Pi_\infty\scr X)$ is $\infty$-connected and identifies local systems with locally constant sheaves on $\scr X$.
\end{proposition}

\begin{proof}
	This is \cite[Theorem A.1.15]{HA}.
\end{proof}

If we were to repeat the proof of Proposition~\ref{prop:lconst} in the world of $(n+1)$-topoi, it would only show that, in a locally $n$-connected $\infty$-topos, $\phi^*$ identifies local systems of $(n-1)$-groupoids (which are always split) with locally constant $(n-1)$-truncated objects. To treat the edge case of local systems of $n$-groupoids, which need not be split, new arguments are needed.

The proof of the following result is the same as the first half of the proof of \cite[Theorem A.1.15]{HA}.

\begin{proposition}\label{prop:splitconst}
	Let $\scr X$ be an $\infty$-topos and $L$ a split local system on $\scr X$. Then the underlying sheaf $\phi^\ast L$ is locally constant. 
\end{proposition}

\begin{proof}
	Let $L$ be represented by $K\in\B(X)$ for some torsor $(X,x)$ in $\scr X$, so that $\phi^\ast L$ is given by the cartesian square
	\begin{tikzmath}
		\diagram{\phi^\ast L & \pi^\ast K \\ \1 & \pi^\ast X\rlap.\\};
		\arrows (11-) edge (-12) (11) edge (21) (12) edge (22) (21-) edge node[above]{$x$} (-22);
	\end{tikzmath}
	Let $\coprod_\alpha U_\alpha\to X$ be an effective epimorphism in $\scr S$ where each $U_\alpha$ is contractible.
	Then $\coprod_\alpha\phi^\ast U_\alpha\to\1$ is an effective epimorphism. There is a commutative diagram in $\Top^\op$
	\begin{tikzmath}
		\diagram{\scr S_{\comma X} & \scr X_{\comma \pi^\ast X} & \scr X_{\comma\1} \\ \scr S_{\comma U_\alpha} & \scr X_{\comma \pi^\ast U_\alpha} & \scr X_{\comma \phi^\ast U_\alpha}\\};
		\arrows (11-) edge (-12) (12-) edge (-13) (11) edge (21) (21-) edge (-22) (22-) edge (-23) (12) edge (22) (13) edge (23);
	\end{tikzmath}
	such that $K$ in the top left corner goes to $\phi^\ast L\times \phi^\ast U_\alpha$ in the bottom right corner.
	Since $\scr S_{\comma U_\alpha}\simeq\scr S$, this shows that $\phi^\ast L\times \phi^\ast U_\alpha$ is constant over $\phi^\ast U_\alpha$. Thus, $\phi^\ast L$ is locally constant.
\end{proof}

\begin{lemma}\label{lem:betterlimits}
	Let $\scr X$ be a locally $n$-connected $\infty$-topos. Then $\pi^*\colon\scr S\to\scr X$ preserves the limits of cofiltered diagrams with $n$-truncated transition maps.
\end{lemma}

\begin{proof}
	Let $K\colon\scr I\to\scr S$ be a cofiltered diagram with $n$-truncated transition maps. Assume without loss of generality that $\scr I$ has a final object $0$. We then have a commutative square
	\begin{tikzmath}
		\diagram{\scr S_{\comma K(0)} & \scr X_{\comma \pi^* K(0)} \\
		\scr S & \scr X \\};
		\arrows (11-) edge node[above]{$\pi^*$} (-12) (11) edge (21) (21-) edge node[above]{$\pi^*$} (-22) (12) edge (22);
	\end{tikzmath}
	where the vertical arrows are the forgetful functors. Since the latter preserve and reflect cofiltered limits, it will suffice to show that $\pi^*\colon (\scr S_{\comma K(0)})_{\leq n}\to (\scr X_{\comma \pi^*K(0)})_{\leq n}$ preserves limits. By descent, this functor can be identified with the functor
	\[\Fun(K(0), \scr S_{\leq n})\to \Fun(K(0),\scr X_{\leq n})\]
	given objectwise by $\pi^*$. This functor preserves limits since they are computed objectwise.
\end{proof}

Recall that a morphism of $\infty$-groupoids $f\colon X\to Y$ is $n$-connected if its fibers are $n$-connected. This is the case if and only if the induced geometric morphism $f_*\colon \B(X)\to \B(Y)$ is $n$-connected. 

\begin{definition}
	Let $-2\leq n\leq \infty$.
	A pro-$\infty$-groupoid is called \emph{$n$-strict} if it can be corepresented by a cofiltered diagram in which the transition maps are $n$-connected.
\end{definition}

\begin{lemma}\label{lem:efftors2}
	Let $\scr X$ be a locally $n$-connected $\infty$-topos.
	Then there exists a coreflective subcategory of $\Tors(\scr X)$ in which all morphisms are $n$-connected.
\end{lemma}

\begin{proof}
	Let $(A,\chi)$ be a torsor. For every morphism of torsors $f\colon (B,\psi)\to(A,\chi)$, consider the unique factorization
	\[B\to e(f)\to A\]
	where $B\to e(f)$ is $n$-connected and $e(f)\to A$ is $n$-truncated. Let $\tilde A$ be the limit of the cofiltered diagram
	\[\Tors(\scr X)_{\comma (A,\chi)}\to \scr S,\quad f\mapsto e(f),\]
	which exists since $\Tors(\scr X)_{\comma (A,\chi)}$ is accessible.
	By construction, this is a diagram with $n$-truncated transition maps. Hence, by Lemma~\ref{lem:betterlimits}, $\pi^*$ preserves the limit of this diagram. In particular, there is an $\tilde A$-torsor $\tilde\chi\colon \1\to \pi^*\tilde A$ which is the limit of the torsors $\1\to \pi^* e(f)$.
	If $\scr C_{(A,\chi)}\subset \Tors(\scr X)_{\comma (A,\chi)}$ denotes the
	image of the localization endofunctor $e$,
	$(\tilde A,\tilde\chi)$ is thus an initial object of $\scr C_{(A,\chi)}$.
	For any $f\colon (B,\psi)\to (A,\chi)$, the functor $\scr C_{(B,\psi)}\to \scr C_{(A,\chi)}$ sending $g$ to $e(f\circ g)$ is left adjoint and hence sends $(\tilde B,\tilde\psi)$ to $(\tilde A,\tilde\chi)$; in particular, the map $\tilde B\to \tilde A$ is $n$-connected. If moreover $B\to A$ is $n$-truncated, then $\tilde B\to \tilde A$ must be an equivalence. Using \cite[Proposition 5.2.7.4]{HTT}, we deduce that $(A,\chi)\mapsto (\tilde A,\tilde \chi)$ is a coreflector with the desired property.
\end{proof}

\begin{samepage}
\begin{proposition}\label{prop:locconn}
	Let $\scr X$ be an $\infty$-topos, let $-2\leq n\leq \infty$, and let $\{X_\alpha\}$ be a family of objects generating $\scr X$ under colimits.
	The following conditions are equivalent:
	\begin{enumerate}
		\item $\scr X$ is locally $n$-connected.
		\item For every $\alpha$, the pro-$n$-groupoid $\tau_{\leq n}\pi_!(X_\alpha)$ is constant.
		\item For every $\alpha$, the pro-$\infty$-groupoid $\pi_! (X_\alpha)$ is $n$-strict.
	\end{enumerate}
\end{proposition}
\end{samepage}

\begin{proof}
	By definition, $\scr X$ is locally $n$-connected if and only if the composition
	\[
	\scr X \xrightarrow{\pi_!}\Pro(\scr S)\xrightarrow{\tau_{\leq n}} \Pro(\scr S_{\leq n})
	\]
	factors through the Yoneda embedding $\scr S_{\leq n}\into \Pro(\scr S_{\leq n})$. Since the latter preserves colimits, we see that $(1)\Leftrightarrow (2)$.
	The implication $(3)\Rightarrow (2)$ is obvious.
	Let us prove $(1)\Rightarrow (3)$. It is clear that $\scr X_{\comma U}$ is locally $n$-connected for every $U\in\scr X$, since $\scr X_{\comma U}\to \scr X$ is étale, so it will suffice to show that $\Pi_\infty\scr X$ is $n$-strict. Let $\scr I\subset\Tors(\scr X)$ be a coreflective subcategory as in Lemma~\ref{lem:efftors2}. By Joyal's criterion \cite[Theorem 4.1.3.1]{HTT}, $\scr I$ is a left cofinal subcategory. Hence, $\Pi_\infty\scr X$ is corepresented by the forgetful functor $\scr I\to\scr S$ and in particular is $n$-strict.
\end{proof}

For $n=-1$, Proposition~\ref{prop:locconn} says that the fundamental pro-$\infty$-groupoid of any $\infty$-topos $\scr X$ can be corepresented by a cofiltered diagram whose transition maps are effective epimorphisms.

\begin{lemma}\label{lem:locnconn}
	If $X\in\Pro(\scr S)$ is $n$-strict, the canonical functor
	\[\Fun(X,\scr S_{\leq n})\to\B(X)_{\leq n}\]
	is fully faithful.
\end{lemma}

\begin{proof}
	Let $X\colon\scr I\to\scr S$ be a corepresentation with $n$-connected transition maps. We claim that each functor
\[\pi_i^*\colon \Fun(X_i,\scr S)\to \B(X)\]
is fully faithful on $n$-truncated objects. Since $\scr I_{\comma i}\to\scr I$ is left cofinal, we may replace $\scr I$ by $\scr I_{\comma i}$ and assume that $i$ is a final object.
Let $p\colon \scr E\to\scr I^\op$ be the topos fibration classified by the functor $\B\circ X\colon\scr I\to\Top$.
By the construction of cofiltered limits of $\infty$-topoi, $\B(X)$ is the $\infty$-category of cartesian sections of $p$. Let $\T(X)$ be the $\infty$-topos of \emph{all} sections of $p$, and factor $\pi_i$ as $\rho_i\circ \sigma$ where $\sigma_*\colon\B(X)\to\T(X)$ is the inclusion and $\rho_{i*}\colon\T(X)\to\Fun(X_i,\scr S)$ is evaluation at $i$. 
Since $\rho_i^*$ is fully faithful, it will suffice to show that, for $F\in\Fun(X_i,\scr S_{\leq n})$, the section $\rho_i^*F$ is already cartesian.
Note that, since $i\in\scr I$ is final, we have $(\rho_i^*F)_j=f^*F$ where $f\colon j\to i$ is the unique map.
Given $f\colon j\to i$, $g\colon k\to j$, and $x\in X_j$, the component at $x$ of the natural transformation $(\rho_i^*F)_j\to g_*(\rho_i^*F)_k$ in $\Fun(X_j,\scr S)$ is the canonical map
\[
F(fx) \to \lim_{y\in g^{-1}(x)} F(fgy) \simeq \Map(g^{-1}(x),F(fx)).
\]
This is an equivalence since $g^{-1}(x)$ is $n$-connected and $F(fx)$ is $n$-truncated.
\end{proof}

\begin{lemma}\label{lem:descent}
	Let $\scr X$ be an $\infty$-topos and let $E\colon \scr I\times\scr J\to \scr X$ be a diagram where $\scr I$ has a final object $e$. Suppose that, for every $i\in\scr I$ and every $j\to k$ in $\scr J$, the square
	\begin{tikzmath}
		\diagram{E(i,j) & E(i,k) \\ E(e,j) & E(e,k) \\};
		\arrows (11-) edge (-12) (21-) edge (-22) (11) edge (21) (12) edge (22);
	\end{tikzmath}
	is cartesian. Then the canonical map
	\[
	\colim_j\lim_i E(i,j)\to \lim_i\colim_j E(i,j)
	\]
	is an equivalence.
\end{lemma}

\begin{proof}
	By descent, the square
	\begin{tikzmath}
		\diagram{E(i,j) & \colim_jE(i,j) \\ E(e,j) & \colim_jE(e,j) \\};
		\arrows (11-) edge (-12) (21-) edge (-22) (11) edge (21) (12) edge (22);
	\end{tikzmath}
	is cartesian for every $i\in\scr I$ and $j\in\scr J$. 
	 Taking the limit over $i$, we obtain a cartesian square
	\begin{tikzmath}
		\diagram{\lim_iE(i,j) & \lim_i\colim_jE(i,j) \\ E(e,j) & \colim_jE(e,j)\rlap. \\};
		\arrows (11-) edge (-12) (21-) edge (-22) (11) edge (21) (12) edge (22);
	\end{tikzmath}
	By universality of colimits, we get the desired equivalence by taking the colimit over $j$ of the left column.
\end{proof}

\begin{samepage}
\begin{lemma}\label{lem:tricky}
	Let $\scr X$ be an $\infty$-topos, $X$ a pro-$\infty$-groupoid, $\phi\colon \scr X\to \B(X)$ a geometric morphism, and $n\geq -2$ an integer. Suppose that:
	\begin{enumerate}
		\item $X$ is $n$-strict;
		\item the composition \[\Fun(X,\scr S_{\leq n})\to \B(X)_{\leq n}\xrightarrow{\phi^*}\scr X_{\leq n}\] is fully faithful.
	\end{enumerate}
	Then $\phi$ is $n$-connected.
\end{lemma}
\end{samepage}

\begin{proof}
	By (1), $X$ is corepresented by a cofiltered diagram $X\colon\scr I\to\scr S$ with $n$-connected transition maps.
	To simplify the notation, we implicitly work with categories of $n$-truncated objects throughout the proof.
	Let $\pi_i\colon\B(X)\to\B(X_i)$ be the canonical projection, and for $j\to i$ in $\scr I$, let $\pi_{ji}\colon \B(X_j)\to \B(X_i)$ be the induced geometric morphism. By Lemma~\ref{lem:locnconn}, $\pi_i$ is $n$-connected, \ie, $\pi_i^*$ is fully faithful. Note that $\pi_{ji}$ is étale and in particular $\pi_{ji}^*$ has a left adjoint $\pi_{ji!}$. Using these facts, we easily verify that $\pi_i^*$ has a left adjoint $\pi_{i!}$ given by
	\begin{equation}\label{eqn:quasietale}
		\pi_{i!}\simeq \colim_j \pi_{ji!}\pi_{j*}.
	\end{equation}
	
	We can assume that $\scr I$ has a final object $e$. Let $E\colon \scr I\times\scr I^\op\to \Fun(\B(X),\B(X))$ be the functor given by
	\[E(i,j)=\pi_i^*\pi_{i!}\pi_j^*\pi_{j*}.\]
	Note that $E(i,j)\to E(j,j)$ is an equivalence for any $i\to j$ in $\scr I$.
	We will show at the end that $E$ satisfies the assumption of Lemma~\ref{lem:descent} (when evaluated at any object in $\B(X)$). It follows that the canonical map
	\[\colim_j E(j,j)\to\lim_i\colim_j E(i,j)\]
	is an equivalence. But the left-hand side is canonically equivalent to the identity functor of $\B(X)$, by \cite[Lemma 6.3.3.6]{HTT}. In other words, for every $L\in \B(X)$,
	\[L\simeq\lim_i \pi_{i}^*\pi_{i!}L.\]
	Applying Lemma~\ref{lem:descent} to $\phi^*\circ E$, we similarly deduce that
	\[\phi^*(L)\simeq\lim_i \phi_i^*\pi_{i!}L,\]
	where $\phi_i=\pi_i\circ\phi$. By assumption (2), $\phi_i^*$ is fully faithful. Therefore,
	\[\phi_{j*}\phi^*\simeq\lim_i\phi_{j*}\phi_i^*\pi_{i!}\simeq \lim_i \pi_{j*}\pi_i^*\pi_{i!}\simeq \pi_{j*}.\]
	 This shows that $\phi^*$ is fully faithful, as desired.
	 
	 We now come back to the claim that, for every $i\in\scr I$ and $k\to j$ in $\scr I$, the square
	 \begin{tikzmath}
	 	\diagram{\pi_i^*\pi_{i!}\pi_j^*\pi_{j*} & \pi_i^*\pi_{i!}\pi_k^*\pi_{k*} \\
		\pi_e^*\pi_{e!}\pi_j^*\pi_{j*} & \pi_e^*\pi_{e!}\pi_k^*\pi_{k*} \\};
		\arrows (11-) edge (-12) (12) edge (22) (11) edge (21) (21-) edge (-22);
	 \end{tikzmath}
	 is cartesian. We see that it suffices to show that the square
	 \begin{tikzmath}
	 	\diagram{\pi_{kj}^*\pi_{kj*} & \id_{\B(X_k)} \\
		\pi_{ke}^*\pi_{je!}\pi_{kj*} & \pi_{ke}^*\pi_{je!}\pi_{kj!} \\};
		\arrows (11-) edge (-12) (12) edge (22) (11) edge (21) (21-) edge (-22);
	 \end{tikzmath}
	 is cartesian, using successively the following facts: $\pi_i^*$ preserves pullbacks, we have the projection formula $A\times_B\pi_{i!}C\simeq \pi_{i!}(\pi_i^*A\times_{\pi_i^*B}C)$ (this follows at once from~\eqref{eqn:quasietale}), $\pi_k^*$ preserves pullbacks, and $\pi_j^*$ and $\pi_k^*$ are fully faithful. Taking the fiber of this square over a point in $X_k$ and using descent in $\scr S$, we are reduced to proving the following statement: if $K$ is a pointed Eilenberg–Mac Lane space of degree $n+1$ and $F\colon K\to\scr S_{\leq n}$ is a functor, then the square
	 \begin{tikzmath}
	 	\diagram{\lim F & F(*) \\
		\lim F & \tau_{\leq n}\colim F \\};
		\arrows (11-) edge (-12) (12) edge (22) (11) edge[-,vshift=1pt] (21) edge[-,vshift=-1pt] (21) (21-) edge (-22);
	 \end{tikzmath}
	 is cartesian. This is obvious if $n\leq 0$. Let us therefore assume that $n\geq 1$, so that $K=K(A,n+1)$ for some abelian group $A$. Let $x$ be a point in $F(*)$. By the long exact sequence associated with the fiber sequence $F(*) \to \colim F \to K$, $F(*)\to\colim F$ is a $\tau_{\leq n-1}$-equivalence and we have an exact sequence
	 \[
	 0 \to \pi_{n+1}(\colim F,x) \to A \to \pi_n(F(*),x) \to \pi_n(\colim F,x)\to 0.
	 \]
	 We must therefore show that, if $x$ is the image of some $y\in \lim F$, then $\pi_{n+1}(\colim F,x) \to A$ is surjective. In fact, $y$ determines a map $K\to \lim F\times K \to \colim F$ which is a section of $(\colim F,x)\to (K,*)$.
\end{proof}

\begin{theorem}
\label{thm:galois}
	Let $-2\leq n\leq\infty$
		and let $\scr X$ be a locally $n$-connected $\infty$-topos. 
		Then the functors
		\begin{equation}
			\label{eqn:phi}
		\Fun(\Pi_\infty\scr X,\scr S_{\leq n}) \to \B(\Pi_\infty\scr X)_{\leq n} \xrightarrow{\phi^*} \scr X_{\leq n}
		\end{equation}
		are fully faithful. They identify $\Fun(\Pi_\infty\scr X,\scr S_{\leq n})$ with the full subcategory of locally constant $n$-truncated sheaves on $\scr X$, and $\B(\Pi_\infty\scr X)_{\leq n}$ with the full subcategory of filtered colimits of locally constant $n$-truncated sheaves on $\scr X$.
		In particular, the geometric morphism $\phi\colon\scr X\to\B(\Pi_\infty\scr X)$ is $n$-connected.
\end{theorem}

\begin{proof}
	The case $n=\infty$ is Proposition~\ref{prop:lconst} (but the following proof also works for $n=\infty$).
	We already know from Proposition~\ref{prop:splitconst} that $\phi^\ast$ sends split local systems to locally constant objects.
	By Proposition~\ref{prop:locconn}, $\Pi_\infty\scr X$ is $n$-strict. By Lemma~\ref{lem:tricky}, it therefore remains to prove that the composition~\eqref{eqn:phi}
	is fully faithful and that every locally constant object is in its image. 
	The description of the image of $\phi^*$ will then follow from the fact that every local system is a filtered colimit of split local systems.

	First we establish a crucial preliminary result.
	By the assumption of local $n$-connectedness, the functor $\pi^\ast\colon\scr S_{\leq n}\to\scr X_{\leq n}$ preserves limits. It follows that for $n$-groupoids $X$ and $Y$,
	\[\pi^\ast\Map(X,Y)\simeq \pi^\ast(\lim_X Y)\simeq \lim_X(\pi^\ast Y)\simeq\Hom(\pi^\ast X,\pi^\ast Y),\]
	where $\Hom$ is the internal mapping object of $\scr X$. Taking global sections, we get, for $n$-groupoids $X$ and $Y$,
	\begin{equation}\label{eqn:cartclosed}
		\Map_{\Pro(\scr S)}(\Pi_\infty\scr X,\Map(X,Y))\simeq\Map_{\scr X}(\pi^\ast X,\pi^\ast Y).
	\end{equation}
	
Let $\scr O\colon \scr X^\op\to\Cat$ be the functor $U\mapsto (\scr X_{\comma U})_{\leq n}$, and let $\scr L\colon \scr X^\op\to\Cat$ be the functor $U\mapsto \Fun(\pi_!U,\scr S_{\leq n})$.
The morphism $\phi$, being natural in $\scr X$, gives rise to a natural transformation $\phi^*\colon \scr L\to\scr O$ whose components are the instances of~\eqref{eqn:phi} for all étale $\scr X$-topoi.
Both $\scr O$ and $\scr L$ are sheaves on $\scr X$, \ie, limit-preserving functors, the former by descent and the latter because $\pi_!$ is left adjoint.
It will therefore suffice to prove the following, for every $U\in\scr X$:
\begin{enumerate}
	\item[(a)] For every $L,M\in\scr L(U)$, there exists a diagram $P\colon\scr A\to\scr X$ with colimit $U$ such that \[\Map_{\scr L(P_\alpha)}(L|P_\alpha, M|P_\alpha) \to \Map_{\scr O(P_\alpha)}(\phi^*L|P_\alpha, \phi^*M|P_\alpha)\] is an equivalence, for all $\alpha\in\scr A$.
	\item[(b)] For every $X\in\scr O(U)$ locally constant, there exists an effective epimorphism $\coprod_\alpha U_\alpha\to U$ such that $X|U_\alpha$ belongs to the essential image of $\scr L(U_\alpha)\to \scr O(U_\alpha)$, for all $\alpha$.
\end{enumerate}
Assertion (b) is obvious since any constant sheaf belongs to the image of $\phi^*$. To prove (a), we choose a torsor $(A,\chi)\in\Tors(\scr X_{\comma U})$ such that $L$ and $M$ are represented by functors $A\to\scr S_{\leq n}$. Let $P\colon A\to\scr X_{\comma U}$ be the associated principal action, whose colimit is $U$. Given $\alpha\in A$, the cartesian square
\begin{tikzmath}
	\diagram{
	P_\alpha & \pi^*(*) \\ U & \pi^*(A) \\
	};
	\arrows (11-) edge (-12) (11) edge (21) (21-) edge node[above]{$\chi$} (-22) (12) edge node[right]{$\pi^*(\alpha)$} (22);
\end{tikzmath}
shows that $L|P_\alpha$ and $M|P_\alpha$ are constant functors $\pi_!P_\alpha\to\scr S_{\leq n}$.
Replacing $\scr X$ by $\scr X_{\comma P_\alpha}$, which is also locally $n$-connected, we are reduced to proving the following statement: for every $X,Y\in\scr S_{\leq n}$, the map
\[
\Map_{\Fun(\pi_!\1,\scr S)}(X,Y)\to \Map_{\scr X}(\pi^*X,\pi^*Y)
\]
is an equivalence. This is exactly our preliminary result~\eqref{eqn:cartclosed}.
\end{proof}

\begin{remark}
	An example of a $1$-topos for which $\phi$ is not connected was constructed by Kennison in \cite[Example 4.16]{Kennison}, showing that the local connectedness assumption in Theorem~\ref{thm:galois} cannot be dropped.
	On the other hand, $\phi$ may be $n$-connected even if $\Pi_n\scr X$ is not constant: if $\scr X=\Shv(\Q)$, then $\phi$ is an equivalence of $\infty$-topoi (but \eqref{eqn:phi} fails to be fully faithful for $n=0$).
\end{remark}

If $\scr X$ is an $\infty$-topos and $n\geq 0$, denote by $\EM_n(\scr X)\subset\scr X_{*/}$ the category of Eilenberg--Mac Lane objects of degree $n$ in $\scr X$.

\begin{corollary}\label{cor:galois}
	Let $\scr X$ be a locally connected $\infty$-topos and let $n\geq 0$. Then the functors
	\[\Fun(\Pi_\infty\scr X,\EM_n(\scr S))\to \EM_n(\B(\Pi_\infty\scr X))\xrightarrow{\phi^*}\EM_n(\scr X)\]
	are fully faithful and identify split local systems (resp.\ local systems) with locally constant sheaves (resp.\ with filtered colimits of locally constant sheaves). 
	If moreover $\scr X$ is locally simply connected, then all local systems of Eilenberg--Mac Lane objects are split.
\end{corollary}

\begin{proof}
	If $n\leq 1$, replace abelian groups by groups or by pointed sets in the following argument. We have a commutative diagram
	\begin{tikzmath}
		\def\colsep{1.5em}
		\diagram{\Fun(\Pi_\infty\scr X,\EM_n(\scr S)) & \EM_n(\B(\Pi_\infty\scr X)) & \EM_n(\scr X) \\ \Fun(\Pi_\infty\scr X,\Ab) & \Ab(\B(\Pi_\infty\scr X)_{\leq 0}) & \Ab(\scr X_{\leq 0}) \\};
		\arrows (11-) edge (-12) (21-) edge (-22) (11) edge node[left]{$\pi_n$} (21) (12-) edge node[above]{$\phi^\ast$} (-13) (22-) edge node[above]{$\phi^\ast$} (-23) (12) edge node[left]{$\pi_n$} (22) (13) edge node[right]{$\pi_n$} (23);
	\end{tikzmath}
	in which the vertical functors are equivalences \cite[Proposition 7.2.2.12]{HTT}. By Theorem~\ref{thm:galois}, the lower row has the desired properties. 
\end{proof}

\begin{definition}
Let $\scr X$ be an $\infty$-topos.
A torsor $(A,\chi)$ in $\scr X$ is called \emph{Galois} if the associated principal action $P\colon A\to \scr X$ is a fully faithful functor.
\end{definition}

 In other words, a Galois torsor in $\scr X$ is a full sub-$\infty$-groupoid of $\scr X$ whose colimit is a final object. For example, if $A\in\scr S$, the universal $A$-torsor is Galois since it corresponds to the Yoneda embedding $A\into \B A$. 
		 If $\scr X$ is the étale topos of a field $k$, then a Galois torsor in $\scr X$ is precisely a finite Galois extension of $k$. In this case, the following corollary shows that the absolute Galois group of $k$ computes $\Pi_1\scr X$.
		 
\begin{corollary}\label{cor:galoisclosure}
	Let $\scr X$ be a locally $n$-connected $\infty$-topos and let $\scr G\subset \Tors_{\leq n+1}(\scr X)$ be the full subcategory spanned by the Galois torsors. Then the inclusion $\scr G\subset \Tors_{\leq n+1}(\scr X)$ is left cofinal.
\end{corollary}

\begin{proof}
	By Lemma~\ref{lem:efftors2}, there exists a left cofinal subcategory $\scr I\subset\Tors_{\leq n+1}(\scr X)$ in which all morphisms are $n$-connected. 
	It will suffice to show that $\scr I\subset \scr G$. If $(A,\chi)$ is a torsor in $\scr I$, the associated principal action is the composition
	\[A\into \B (A)_{\leq n}\to \B(\Pi_{n+1}\scr X)_{\leq n}\xrightarrow{\phi^*} \scr X_{\leq n}.\]
	The second arrow is fully faithful by Lemma~\ref{lem:locnconn} and the third by Theorem~\ref{thm:galois}.
\end{proof}

An $\infty$-topos $\scr X$ is called \emph{$n$-Galois} if it is $n$-localic, locally $(n-1)$-connected, and generated under colimits by the images of its Galois torsors $A\into\scr X$ with $A\in\scr S_{\leq n}$.
For $n=1$, this recovers the usual notion of Galois topos \cite[\S3]{Moerdijk}, except that we do not require connectedness.
We then have the following generalization of \cite[Theorem 3.2]{Moerdijk}:

\begin{samepage}
\begin{corollary}
	Let $\scr X$ be an $n$-localic $\infty$-topos, $-1\leq n\leq \infty$. The following are equivalent:
	\begin{enumerate}
		\item $\scr X$ is $n$-Galois;
		\item $\scr X$ is locally $(n-1)$-connected and generated under colimits by its locally constant objects;
		\item $\scr X$ is locally $(n-1)$-connected and generated under colimits and finite limits by its locally constant objects;
		\item $\scr X$ is locally $(n-1)$-connected and $\phi\colon\scr X\to\B(\Pi_{n}\scr X)$ is an equivalence of $\infty$-topoi;
		\item There exist an $(n-1)$-strict pro-$n$-groupoid $X$ and an equivalence $\scr X\simeq\B(X)$.
	\end{enumerate}
\end{corollary}
\end{samepage}

\begin{proof}
	(1) $\Rightarrow$ (2). It suffices to note that, if $P\colon A\to \scr X$ is a torsor, then $P_\alpha$ is locally constant for all $\alpha\in A$. (2) $\Rightarrow$ (3). Obvious. (3) $\Rightarrow$ (4). Follows from Theorem~\ref{thm:galois}. (4) $\Rightarrow$ (5). Follows from Proposition~\ref{prop:locconn}. (5) $\Rightarrow$ (1). Let $X\colon\scr I\to\scr S_{\leq n}$ be a cofiltered diagram with $(n-1)$-connected transition maps. It is clear that $\B(X)$ is generated under colimits by the images of the torsors
	\[X_i\into \B(X_i)\to \B(X).\]
	These are Galois torsors by Lemma~\ref{lem:locnconn}. That $\B(X)$ is locally $(n-1)$-connected was verified at the beginning of the proof of Lemma~\ref{lem:tricky}.
\end{proof}

\section{Finite Galois theory} 
\label{sec:finitegalois}

In this short section, we give a classification of \emph{finite} locally constant sheaves in an $\infty$-topos.
Compared to the classification of $n$-truncated locally constant sheaves established in \S\ref{sec:galois}, the main difference is that we do not need any assumption of local connectedness.

Recall from \cite[\S2.3]{DAG13} that an $\infty$-groupoid is \emph{finitely constructible} if it is truncated and coherent, \ie, has finite homotopy groups. We denote by $\scr S^\fc\subset\scr S$ the full subcategory of finitely constructible $\infty$-groupoids.

\begin{definition}
	Let $\scr X$ be an $\infty$-topos. An object $X\in\scr X$ is \emph{finite locally constant} if there exists an effective epimorphism $\coprod_\alpha U_\alpha\to\1$ and finitely constructible $\infty$-groupoids $X_\alpha\in\scr S^\fc$ such that $X\times U_\alpha\simeq \pi^*X_\alpha\times U_\alpha$ in $\scr X_{\comma U_\alpha}$, for all $\alpha$.
\end{definition}

Let $\scr S_{<\infty}\subset\scr S$ be the full subcategory of truncated $\infty$-groupoids. Since it is stable under finite limits, we have an adjunction
\begin{tikzmath}
	\diagram{\Pro(\scr S) & \Pro(\scr S_{<\infty})\rlap. \\};
	\arrows (11-) edge[vshift=\dbl] node[above]{$\tau_{<\infty}$} (-12) (-12) edge[c->,vshift=\dbl] (11-);
\end{tikzmath}

\begin{lemma}\label{lem:grothfc}
	Let $\scr X$ be an $\infty$-topos and let $A\in\scr S$ be coherent. Then the canonical map
	\[
	\Pi_\infty(\scr X_{\comma \pi^*A}) \to A\times \Pi_\infty\scr X
	\]
	is a $\tau_{<\infty}$-equivalence.
\end{lemma}

\begin{proof}
	By descent, we have
	\[
	\Pi_\infty(\scr X_{\comma \pi^*A}) \simeq \colim_{a\in A} \Pi_\infty\scr X.
	\]
	To conclude, we observe that for $X\in\Pro(\scr S)$ and $A\in\scr S$ coherent, the canonical map $\colim_{A}X \to A\times X$ is a $\tau_{<\infty}$-equivalence. Indeed, by \cite[Proposition 2.3.9]{DAG13}, $\Map(A,\ph)\colon\scr S_{\leq n}\to\scr S$ preserves filtered colimits for all finite $n$. If $Y\in\scr S_{<\infty}$, we therefore have
	\begin{multline*}
	\Map(\colim_AX,Y)\simeq \Map(A,\Map(X,Y))\simeq \Map(A,\colim_{i}\Map(X_i,Y))
	\\
	\simeq\colim_i\Map(A,\Map(X_i,Y))\simeq \colim_i\Map(A\times X_i,Y)\simeq \Map(A\times X,Y),
	\end{multline*}
	which proves the result.
\end{proof}

\begin{theorem}\label{thm:finitegalois}
	Let $\scr X$ be an $\infty$-topos. Then the geometric morphism $\phi\colon\scr X\to \B(\Pi_\infty\scr X)$ induces a fully faithful functor \[\Fun(\Pi_\infty\scr X,\scr S^\fc)\into \scr X\] whose essential image is the subcategory of finite locally constant sheaves on $\scr X$.
\end{theorem}

\begin{proof}
	By the proof of Proposition \ref{prop:splitconst}, this functor lands in the subcategory of finite locally constant sheaves.
	Since the functors $U\mapsto \scr X_{\comma U}$ and $U\mapsto \Fun(\pi_!U,\scr S^\fc)$ are sheaves on $\scr X$,
	we are reduced as in the proof of Theorem~\ref{thm:galois} to proving the following statement: for every $A,B\in\scr S^\fc$, the map
	\[
	\Map_{\Fun(\pi_!\1,\scr S)}(A,B) \to \Map_{\scr X}(\pi^*A,\pi^*B)
	\]
	is an equivalence. 
	Replacing $\scr X$ by $\scr X_{\comma \pi^*A}$ and using Lemma~\ref{lem:grothfc}, we can assume that $A=*$. The result then follows from Proposition~\ref{prop:observation}.	
\end{proof}

\section{The étale homotopy type of Artin–Mazur–Friedlander} 
\label{sec:examples}

Let $\scr C$ be a small Grothendieck site with a final object. We denote by $\Shv(\scr C)$ the $\infty$-topos of sheaves of $\infty$-groupoids on $\scr C$ and by $\Shv(\scr C)^\wedge$ its hypercompletion.
Let $\HC(\scr C)$ be the category of hypercovers of $\scr C$, in the sense of \cite[Definition 4.2]{DHI}. This is a full subcategory of the category $\Fun(\scr C^\op,\Set_\Delta)$ of simplicial presheaves on $\scr C$, and as such it is enriched in simplicial sets. We let $\pi\HC(\scr C)$ be the category obtained from $\HC(\scr C)$ by identifying simplicially homotopic morphisms. 

The constant simplicial presheaf functor $\Set_\Delta\to \Fun(\scr C^\op,\Set_\Delta)$ has a (simplicially enriched) left adjoint, and we denote by
\[
\Pi\scr C\colon \HC(\scr C)\to \Set_{\Delta}
\]
its restriction to $\HC(\scr C)$.
We would like to regard $\Pi\scr C$ as a cofiltered diagram of $\infty$-groupoids.
The issue that usually arises at this point is that the category $\HC(\scr C)$ is not cofiltered, and while $\pi \HC(\scr C)$ is cofiltered, $\Pi\scr C$ does not identify simplicially homotopic morphisms. This led Artin and Mazur \cite[\S9]{AM} to consider instead the induced functor $\pi\HC(\scr C)\to \h\scr S$, which is a cofiltered diagram in the homotopy category $\h\scr S$.
However, $\HC(\scr C)$ is cofiltered as a \emph{simplicially enriched} category, in the following sense:
\begin{itemize}
	\item Any finite collection of hypercovers admits a common refinement.
	\item For any inclusion of finite simplicial sets $K\subset L$ and any map $K\to\Map_{\Delta}(V,U)$, there exists a refinement $W\to V$ of $V$ such that the induced map $K\to\Map_{\Delta}(W,U)$ extends to $L$.
\end{itemize}
Both statements follow from \cite[Proposition 5.1]{DHI}.
They imply that the $\infty$-category associated with $\HC(\scr C)$ is cofiltered (by \cite[Proposition 5.3.1.13]{HTT} and an easy argument using Kan's $\Ex^\infty$ functor).

\begin{proposition}\label{prop:verdier}
	Let $\scr C$ be a small site with a final object.
	 Then $\Pi_\infty\Shv(\scr C)^\wedge$ is corepresented by the simplicially enriched cofiltered diagram
	 $\Pi\scr C\colon \HC(\scr C)\to \Set_\Delta$.
\end{proposition}

\begin{proof}
	Let $K$ be a Kan complex.
	We must define an equivalence of $\infty$-groupoids
	\[\hocolim_{U\in\HC(\scr C)}\Map_\Delta(\Pi\scr C(U),K)\simeq \Map(\1,\pi^*K),\]
	 natural in $K$, where the right-hand side is a mapping space in $\Shv(\scr C)^\wedge$.
	 Let $\underline K$ be the constant simplicial presheaf on $\scr C$ with value $K$. By definition of $\Pi\scr C$, we have a natural isomorphism of Kan complexes $\Map_\Delta(\Pi\scr C(U),K)\simeq \Map_\Delta(U,\underline K)$.
	The localization functor $\Fun(\scr C^\op,\Set_\Delta)\to \Shv(\scr C)^\wedge$ sends each hypercover to a final object and hence induces a natural map of $\infty$-groupoids
	\begin{equation}\label{eqn:Verdier}
	\hocolim_{U\in\HC(\scr C)}\Map_\Delta(U,\underline K)\to \Map(\1,\pi^*K).
	\end{equation}
	Let $x\colon V\to \underline K$ be a point in the left-hand side. Since $\HC(\scr C)$ is cofiltered and $\pi_n\colon\scr S_{*\comma}\to\scr S_{\leq 0}$ preserves filtered colimits, we have an isomorphism
	\[
	\pi_n(\hocolim_{U\in\HC(\scr C)}\Map_\Delta(U,\underline K),x) \simeq \colim_{U\in\pi\HC(\scr C)_{\comma V}} \pi_n(\Map_\Delta(U,\underline K), x).
	\]
	Applying the generalized Verdier hypercovering theorem \cite[Theorem 8.6]{DHI} to the constant simplicial presheaf $\underline K$, we therefore deduce that~\eqref{eqn:Verdier} is an equivalence, as desired.
\end{proof}

\begin{remark}
	The pro-$\infty$-groupoid $\Pi_\infty\Shv(\scr C)$ does not admit such an explicit model in general.
	Note however that, for any $\infty$-topos $\scr X$, the map $\Pi_{\infty}\scr X^\wedge\to\Pi_\infty\scr X$ is a $\tau_{<\infty}$-equivalence, since truncated objects are hypercomplete. 
	In particular, $\Pi_n\scr X^\wedge\simeq \Pi_n\scr X$ for any finite $n$.
\end{remark}

\begin{lemma}\label{lem:1cofinal}
	Let $\scr I$ and $\scr J$ be cofiltered $\infty$-categories and let $f\colon \scr J\to \scr I$ be a functor. Then $f$ is left cofinal if and only if $\h f\colon \h\scr J\to\h\scr I$ is left $1$-cofinal.
\end{lemma}

\begin{proof}
	Assume that $\h f$ is left $1$-cofinal, and let $p\colon \scr I^\op\to\scr S$ be a diagram. We must show that the map
	\[\colim_{j\in\scr J}pf(j)\to\colim_{i\in\scr I}p(i)\] 
	is an equivalence, \ie, induces isomorphisms on homotopy groups. This follows from the assumption and the fact that homotopy groups $\pi_n\colon\scr S_{*\comma}\to\scr S_{\leq 0}$ preserve filtered colimits. The other implication is obvious.
\end{proof}

\begin{corollary}
	\label{cor:verdier}
	Let $\scr C$ be a small site with a final object, $\scr I$ a cofiltered category, and $f\colon \scr I\to \HC(\scr C)$ a functor such that the composition
	\[
	\scr I\stackrel f\to\HC(\scr C)\to \pi\HC(\scr C)
	\]
	is left $1$-cofinal. Then the cofiltered diagram $\Pi\scr C\circ f$ corepresents $\Pi_\infty\Shv(\scr C)^\wedge$.
\end{corollary}

\begin{proof}
	By Lemma~\ref{lem:1cofinal}, $f$ induces a left cofinal functor of $\infty$-categories. The result now follows from Proposition~\ref{prop:verdier}.
\end{proof}

We now turn to the étale site $\ET_X$ of a scheme $X$, \ie, the category of étale $X$-schemes equipped with the étale topology.
 We denote by $X_\et$ the $\infty$-topos of sheaves of $\infty$-groupoids on $\ET_X$ and by $X_\et^\wedge$ its hypercompletion (although $\ET_X$ is not small, it is easy to show that $X_\et$ is indeed an $\infty$-topos).
The abelian cohomology of either $\infty$-topos is the étale cohomology of $X$.

If $X_\Zar$ denotes the $\infty$-topos of sheaves on the small Zariski site of $X$, there is an obvious geometric morphism $X_\et\to X_\Zar$ which induces an equivalence on $(-1)$-truncated objects. In particular, $\Pi_0(X_\et)\simeq \Pi_0(X_\Zar)$.
By Proposition~\ref{prop:locconn}, we deduce that $X_\et$ is locally connected if and only if the underlying topological space of $X$ is locally connected (as this property persists after étale extension).
Note however that $X_\et$ is rarely locally simply connected, since the étale fundamental pro-groupoid is often nonconstant. 
If $X$ is locally connected, then:
\begin{itemize}
	\item Locally constant sheaves of sets on $\ET_X$ are classified by functors $\Pi_1 X_\et\to \scr S_{\leq 0}$ (Theorem~\ref{thm:galois}).
	\item If $A$ is a locally constant sheaf of abelian groups on $\ET_X$ classified by $L\colon\Pi_1 X_\et\to \Ab$, then $H^*(X_\et,A)\simeq H^*(\Pi_\infty X_\et,L)$ (apply Proposition~\ref{prop:observation} to the Eilenberg–Mac Lane objects $K(L,n)$).
\end{itemize}
On the other hand, with no assumptions on the scheme $X$, we have:
\begin{itemize}
	\item Finite locally constant sheaves on $\ET_X$ are classified by functors $\Pi_\infty X_\et \to \scr S^\fc$ (Theorem~\ref{thm:finitegalois}).
\end{itemize}

In \cite[\S4]{Friedlander:1982}, Friedlander defines the \emph{étale topological type} of a locally connected scheme $X$:
it is a pro-simplicial set given by a diagram
\[
\scr I\to \HC'(\ET_X)\subset \HC(\ET_X)\xrightarrow{\Pi\ET_X} \Set_\Delta,
\]
where $\scr I$ is some small cofiltered poset (whose definition is specific to the étale site) and $\HC'(\ET_X)\subset \HC(\ET_X)$ is the full subcategory spanned by the \emph{basal} hypercovers \cite[\S9]{DHI} that are moreover degreewise coproducts of \emph{connected} étale $X$-schemes. By \cite[Theorem 9.6]{DHI}, any hypercover of $\ET_X$ is refined by a basal hypercover, and since $X$ is locally connected, it is clear that any such is further refined by a hypercover in $\HC'(\ET_X)$. In particular, the inclusion $\pi\HC'(\ET_X)\subset\pi\HC(\ET_X)$ is left $1$-cofinal.
The key property of Friedlander's construction is that the composite functor
\[\scr I\to\HC'(\ET_X)\to\pi\HC'(\ET_X)\]
is left $1$-cofinal \cite[p.~38]{Friedlander:1982}.
Combining these facts and applying Corollary~\ref{cor:verdier}, we deduce the following:

\begin{corollary}
	Let $X$ be a locally connected scheme. Then the étale topological type of $X$, as defined by Friedlander,
	corepresents the shape of $X_\et^\wedge$.
\end{corollary}

\providecommand{\bysame}{\leavevmode\hbox to3em{\hrulefill}\thinspace}


\begin{thebibliography}{AGV72}

\bibitem[AGV72]{SGA4-1}
M.~Artin, A.~Grothendieck, and J.-L. Verdier, \emph{Th{\'e}orie des topos et
  cohomologie {\'e}tale des sch{\'e}mas}, Lecture Notes in Mathematics, vol.
  269, Springer, 1972

\bibitem[AM69]{AM}
M.~Artin and B.~Mazur, \emph{{\'E}tale Homotopy}, Lecture Notes in Mathematics,
  vol. 100, Springer, 1969

\bibitem[BD81]{BarrDiaconescu}
M.~Barr and R.~Diaconescu, \emph{On locally simply connected toposes and their
  fundamental groups}, Cahiers Top. G{\'e}o. Diff. Cat{\'e}goriques \textbf{22}
  (1981), pp.~301--314

\bibitem[BK72]{BK}
A.~K. Bousfield and D.~M. Kan, \emph{Homotopy Limits, Completions and
  Localizations}, Lecture Notes in Mathematics, vol. 304, Springer, 1972

\bibitem[Bun92]{Bunge}
M.~Bunge, \emph{Classifying toposes and fundamental localic groupoids},
  Category Theory 1991 (R.~A.~G. Seely, ed.), CMS Conf. Proc., vol.~13, AMS,
  1992

\bibitem[DHI04]{DHI}
D.~Dugger, S.~Hollander, and D.~C. Isaksen, \emph{Hypercovers and simplicial
  presheaves}, Math. Proc. Cambridge Philos. Soc. \textbf{136} (2004), no.~1,
  pp.~9--51, preprint \href{http://arxiv.org/abs/math/0205027}{{\sf
  arXiv:math/0205027 [math.AT]}}

\bibitem[Dub08]{Dubuc}
E.~J. Dubuc, \emph{The fundamental progroupoid of a general topos}, J. Pure
  Appl. Algebra \textbf{212} (2008), no.~11, pp.~2479--2492

\bibitem[Dub10]{Dubuc2}
\bysame, \emph{Spans and simplicial families}, 2010,
  \href{http://arxiv.org/abs/1012.6001}{{\sf arXiv:1012.6001 [math.CT]}}

\bibitem[Fri82]{Friedlander:1982}
E.~M. Friedlander, \emph{{\'E}tale Homotopy of Simplicial Schemes}, Annals of
  Mathematical Studies, vol. 104, Princeton University Press, 1982

\bibitem[Gro75]{GrothendieckBreen}
A.~Grothendieck, Letter to L.~Breen, 1975,
  \url{http://webusers.imj-prg.fr/~leila.schneps/grothendieckcircle/Letters/breen1.html}

\bibitem[Isa07]{Isaksen:2007}
D.~C. Isaksen, \emph{Strict model structures for pro-categories}, Categorical
  decomposition techniques in algebraic topology, Progress in Mathematics, vol.
  215, Birkhauser, 2007, pp.~179--198, preprint
  \href{http://arxiv.org/abs/math/0108189}{{\sf arXiv:math/0108189 [math.AT]}}

\bibitem[Ken91]{Kennison}
J.~F. Kennison, \emph{The fundamental localic groupoid of a topos}, J. Pure
  Appl. Algebra \textbf{77} (1991), no.~1, pp.~67--86

\bibitem[Lur09]{HTT}
J.~Lurie, \emph{Higher Topos Theory}, Annals of Mathematical Studies, vol. 170,
  Princeton University Press, 2009

\bibitem[Lur11]{DAG13}
\bysame, \emph{Derived Algebraic Geometry XIII: Rational and $p$-adic Homotopy
  Theory}, 2011,  \url{http://www.math.harvard.edu/~lurie/papers/DAG-XIII.pdf}

\bibitem[Lur16]{HA}
\bysame, \emph{Higher Algebra}, 2016,
  \url{http://www.math.harvard.edu/~lurie/papers/HA.pdf}

\bibitem[MS71]{MS}
S.~Marde{\v s}i{\'c} and J.~Segal, \emph{Shapes of compact and ANR-systems},
  Fund. Math. \textbf{72} (1971), pp.~41--59

\bibitem[Moe89]{Moerdijk}
I.~Moerdijk, \emph{Prodiscrete groups and Galois toposes}, Indagat. Math.
  (Proceedings) \textbf{92} (1989), no.~2, pp.~219--234

\bibitem[TV03]{TV}
B.~To{\"e}n and G.~Vezzosi, \emph{Segal topoi and stacks over Segal
  categories}, 2003, \href{http://arxiv.org/abs/math/0212330v2}{{\sf
  arXiv:math/0212330v2 [math.AT]}}

\end{thebibliography}
\end{document}